\documentclass[11pt]{article}

\voffset= - 1.2in
\hoffset= - 1.0in

\def\be{\begin{equation}}
\def\ee{\end{equation}}
\def\bea{\begin{eqnarray}}
\def\eea{\end{eqnarray}}

\title{{\bf
On projective group properties of the $6D$ pseudo-Riemannian space
}
\vspace{.5cm}}
\author{{\bf Zolfira Zakirova}\footnote{E-mail:
zolya@itep.ru}
\date{ } \\
{\small {\it Kazan State Energy University, Kazan, Russia}}
}

\begin{document}

\maketitle

\vspace{-6.5cm}

\begin{center}
\hfill ITEP/TH-46/05\\
\end{center}

\vspace{3.5cm}

\begin{abstract}
We study the six-dimensional pseudo-Riemannian spaces with two time-like coordinates
that admit non-homothetic infinitesimal projective transformations.
The metrics are manifestly obtained and the projective group properties are determined.
We also find a generic defining of projective motion in the 6-dimensional rigid
$h$-space.
\end{abstract}
\def\thefootnote{\arabic{footnote}}

\section{Introduction}     %\section*{Introduction}

The problem of defining $2D$ Riemannian manifolds which admit
projective motions, i. e. continuous transformation groups preserving
geodesics dates back to  S.Lie who considered it in \cite{[1]}.\footnote{Another
important result was obtained by A.Z.Petrov \cite{[2]}, who
classified geodesically equivalent pseudo-Riemannian spaces $V^3$.}
In more recent times, that was
A.V.Aminova who has got a complete solution to this problem \cite{[3]}.
For the Riemannian manifolds of dimension greater than 2
the same problem has been solved
by G.Fubini \cite{[4]} and  A.S.Solodovnikov \cite{[5]}.

Note that they essentially used in their studies positive definiteness of metrics under
consideration. Without this positivity condition, considering
pseudo-Riemannian spaces, the problem is much more complicated and requires
absolutely new method of solution.

In later paper, \cite{[6]} A.V.Aminova has classified
all the Lorentzian manifolds
of dimension $\geq 3$ that admit nonhomothetic projective or affine
infinitesimal transformations. In each case, there were determined the
corresponding maximal projective and affine Lie algebras.

\paragraph{The General problem} is to classify the n-dimensional
pseudo-Riemannian spaces admitting
projective motions, i. e. continuous transformation groups preserving
geodesics.

This problem is not solved yet for pseudo-Riemannian spaces with arbitrary
signature.

\paragraph{Our concrete problem} here is to study a projective group properties of a 6-dimensional
pseudo-Riemannian space with signature [+ + - - - -], which admits
projective motions, i. e. continuous transformation groups preserving
geodesics.

\paragraph{The Method} we use here is the method of the skew-normal frames and general approach to
investigation of  projective motions of
pseudo-Riemannian manifolds due to A. V. Aminova.

\section{The metrics of the rigid $h$-spaces}

Let me remind that vector field $X$ on a pseudo-Riemannian manifold $(M,g)$ is an
infinitesimal projective transformation if and only if \cite{[7]}
\be
L_{X} g=h,
\ \ \ \ \ \
\hbox{{\it generalized Killing equation}}
\ee
\be
\nabla h(Y, Z, W)=2g(Y, Z)W \varphi+g(Y, W)Z \varphi+g(Z, W)Y \varphi,
\ee
\centerline{{\it Eisenhart equation}}
where  $(Y, Z, W) \in T(M)$, $\varphi=\frac{1}{n+1} {\rm div} X$.

In order to find a pseudo-Riemannian space admitting a nonhomothetic
infinitesimal projective transformation, one needs to integrate the
Eisenhart equation.

Let us introduce the following definitions.

\bigskip

\noindent
{\bf Definition 1:}  Pseudo-Riemannian manifolds for which there exist nontrivial
solutions
$h\ne cg$ to the Eisenhart equation are called
{\it{$h$-spaces}}.

\bigskip

The problem of determining such spaces depends on the type of
$h$-space, i. e. on the type of the bilinear form
$L_{X}g$ determined by the characteristic of the $\lambda$-matrix
$( h-\lambda g)$. If the Segre characteristics of the tensor $L_{X}g$ is
$[abc...]$, we call the corresponding $h$-space the space of $[abc...]$ type.
The number of possible types depends on the
dimension and the signature of the pseudo-Riemannian space.
We restrict ourselves to considering the rigid $h$-spaces.

\bigskip

\noindent
{\bf Definition 2:}  $h$-spaces with distinct bases of prime divisors of the
$\lambda$-matrix are called {\it{rigid $h$-spaces}}.

\bigskip

So, we investigate the rigid $h$-spaces, i. e. $h$-spaces of the $[111111]$, $[21111]$, $[2211]$, $[3111]$,
$[321]$, $[33]$, $[411]$ and $[51]$ types.

Using the technics integrating in skew-normal frames, developed by A.V.Aminova,
one finds the metrics of these  $h$-spaces.

In particular,
the metric of the $h$-space of $[2211]$ type is
\be
g_{ij}dx^idx^j=e_2 (f_4-f_2)^2  \Pi_{\sigma} (f_{\sigma}-f_2) \lbrace 2A  dx^1dx^2-A^2 \Sigma_1 (dx^2)^2 \rbrace+
\ee
$$
e_4 (f_2-f_4)^2  \Pi_{\sigma} (f_{\sigma}-f_4) \lbrace 2 \tilde{A}  dx^3dx^4-\tilde{A}^2 \Sigma_2 (dx^4)^2 \rbrace+
\sum_{\sigma} e_{\sigma} {\Pi'_i}(f_i-f_{\sigma}) (dx^{\sigma})^2,
$$
where
$$
A=\epsilon x^1+\theta(x^2),
\quad
\tilde A=\tilde{\epsilon} x^3+\omega(x^4),
$$
$$
\Sigma_1=2 (f_4-f_2)^{-1}+ \sum_{\sigma}(f_{\sigma}-f_2)^{-1},
\quad
\Sigma_2=2 (f_2-f_4)^{-1}+ \sum_{\sigma}(f_{\sigma}-f_4)^{-1},
$$
$f_2=\epsilon x^2$, $f_4=\tilde{\epsilon}x^4+a$, $\epsilon, \tilde{\epsilon}=0, 1$,
$a$ is a constant which is nonzero when $\tilde{\epsilon}=0$,
$f_{\sigma}(x^{\sigma})$, $\theta(x^2)$, $\omega(x^4)$ are arbitrary functions, $e_i=\pm 1$, $\sigma=5, 6$.

The tensor $h_{ij}$ of the $h$-space of $[2211]$ type is
\be
h_{ij} dx^i dx^j=2 f_2 g_{12}dx^1dx^2+(f_2 g_{22}+Ag_{12})(dx^2)^2+
2 f_4 g_{34}dx^3dx^4+
\ee
$$
(f_4 g_{44}+\tilde{A}g_{34})(dx^4)^2+\sum_{\sigma} f_{\sigma} g_{\sigma \sigma}
(dx^{\sigma})^2+
(2 f_2+2 f_4+\sum_{\sigma} f_{\sigma}+c)g_{ij},
$$
For every solution $h_{ij}$ of the geodesic Eisenhart equation, there is a
quadratic first  integral
\be
(h_{ij}-4\varphi g_{ij})\dot{x}^i \dot{x}^j={\rm const},
\ee
where $\dot{x}^i$ is the tangent vector to the geodesic.

For the other rigid $h$-spaces we have  similar results.

\section{On projective group properties of the
\\
6-dimensional pseudo-Riemannian space}

To move further, we need a necessary and sufficient
condition of constant curvature of this $h$-spaces.
A necessary and sufficient
of condition of constant curvature is expressed by the formula
\be
R_{jkl}^i=K({\delta_k}^i g_{jl}-{\delta_l}^i g_{jk}),
\quad
K={\rm const}.
\ee
Calculating components of the curvature tensor of the rigid $h$-spaces and
substituting into this equality, one obtains a necessary and sufficient
conditions of constant curvature of this $h$-space.

In particular, for the $h$-space of $[2211]$ type
\be
\rho_p-\rho_{\sigma p}=\rho_p-\rho_{pq}=\epsilon=\tilde \epsilon=0
\quad
(p\ne q, p, q=2, 4, \sigma=5, 6),
\ee
where
$$
\rho_p=-\frac{1}{4} \sum_{\sigma} \frac{(f'_{\sigma})^2}{(f_{\sigma}-f_p)^2
g_{\sigma \sigma}},
\quad
\rho_{pq}=-\frac{1}{4} \sum_{\sigma} \frac{(f'_{\sigma})^2}
{(f_{\sigma}-f_p)(f_{\sigma}-f_q)g_{\sigma \sigma}},
$$
$$
\rho_{\sigma p}=-\frac{1}{4}\frac{(f'_{\sigma})^2}{(f_{\sigma}-f_p)
g_{\sigma \sigma}} \lbrace \frac{2 {f''}_{\sigma}}{(f'_{\sigma})^2}-
\frac{1}{f_{\sigma}-f_p}+
\sum_{i, i \ne \sigma} (f_i-f_{\sigma})^{-1} \rbrace-
$$
$$
\frac{1}{4} \sum_{\gamma, \gamma \ne \sigma} \frac{(f'_{\gamma})^2}
{(f_{\gamma}-f_p)(f_{\gamma}-f_{\sigma})g_{\gamma \gamma}}.
$$

Further, investigating the Eisenhart equations and their integrability
conditions for each obtained rigid $h$-spaces
we prove some theorems,
when give important information about
structure projective Lie algebra in the rigid $h$-spaces.

\bigskip

\noindent
{\bf Theorem 1.} {\it
Any defining function of
projective motion in rigid $h$-spaces
of nonconstant curvature can be presented as $\phi=a_1 \varphi$,
where $a_1$ is a constant.}

\bigskip

\noindent
{\bf Theorem 2.} {\it Any covariantly constant symmetric
tensor $b_{ij}$ in 6-dimensional
rigid $h$-spaces of nonzero curvature is proportional to the
fundamental tensor, i. e. $b_{ij}=a_2 g_{ij}$, where $a_2$ is a constant.}

\bigskip

From this theorems, one immediately obtains

\bigskip

\noindent
{\bf Theorem 3.} {\it The affine group of 6-dimensional
rigid $h$-spaces of non-constant curvature consists of
homothetics.}

\bigskip

These theorems and linearity of the Eisenhart
equation give the general solution of the Eisenhart
equation in the rigid $h$-spaces of non-constant curvature
in the form
\be
a_1 h_{ij}+a_2 g_{ij}
\ee
with two arbitrary constants $a_1, a_2$.

Hence, one obtains

\bigskip

\noindent
{\bf Theorem 4.} {\it All projective motions of a
{\rm 6}-dimensional rigid $h$-spaces of non-constant curvature
are obtained by integrating the equation
\be
L_{\xi} g_{ij}=\xi_{i,j}+\xi_{j,i}=a_1 h_{ij}+a_2g_{ij}.
\ee}

\bigskip

It leads to the following important group characteristic of
rigid $h$-spaces:

\bigskip

\noindent
{\bf Theorem 5.} {\it If rigid $h$-spaces
of non-constant curvature admit a nonhomothetic projective Lie algebra $P_r$,
then this algebra contains the subalgebra $H_{r-1}$ of infinitesimal
homothetics of dimension $r-1$.}

\bigskip

I am grateful to A.V.Aminova for the constant encouragement and
discussions and to C.BurdiK for kind hospitality in Prague and the
very nice Conference. The work is partially supported by RFBR grant 04-02-16538- .


\bigskip

\begin{thebibliography}{9}
\bibitem{[1]} M. G. Konigs: Appl. II to G. Darboux. Lecons sur la theorie
generalle des surfaces. IV. (1896), p. 368.
\bibitem{[2]} A.Z. Petrov: Uch. zap. Kazan. un-ta. (1949), Vol. 109, no. 3, p. 7 (in Russian).
\bibitem{[3]} A. V. Aminova: Iz. VUZ {\bf 30} (1990), no. 6, p. 1 (in Russian).
\bibitem{[4]} G. Fubini: Mem. Acc. Torino. Cl. Fif. Mat. Nat. (1903) Vol. 53, N2, p. 261.
\bibitem{[5]} A. S. Solodovnikov: Uspekhi Mat. Nauk (1956), no. 11, p. 45 (in Russian).
\bibitem{[6]} A. V. Aminova: Uspekhi Mat. Nauk {\bf 50} (1995), no. 1(301), p. 69 (in Russian).
\bibitem{[7]} L. P. Eisenhart: {\it Riemannian geometry} M.:IL (1948) (in Russian).

\end{thebibliography}
\end{document}